\documentclass[11pt, a4paper]{article}
\usepackage{times}
\usepackage{a4wide}
\usepackage[dvips]{hyperref}
\usepackage[british]{babel}
\usepackage{enumerate}
\usepackage{amsmath, amscd, amsfonts, amssymb, latexsym, theorem}
\usepackage{xypic}
\usepackage{wrapfig}
\usepackage{graphicx}
\usepackage[T1]{fontenc}
\usepackage[latin1]{inputenc}

\theoremstyle{change}
{\theorembodyfont{\itshape}   \newtheorem{theorem}{Theorem.}[section]}
{\theorembodyfont{\itshape}   \newtheorem{lemma}[theorem]{Lemma.}}
{\theorembodyfont{\itshape}   \newtheorem{definition}[theorem]{Definition.}}
{\theorembodyfont{\itshape}   \newtheorem{remark}[theorem]{Remark.}}
{\theorembodyfont{\itshape}   \newtheorem{proposition}[theorem]{Proposition.}}
{\theorembodyfont{\itshape}   \newtheorem{corollary}[theorem]{Corollary.}}
{\theorembodyfont{\itshape}   \newtheorem{notation}[theorem]{Notation.}}

\newcommand{\CC}{\mathbb{C}}
\newcommand{\QQ}{\mathbb{Q}}
\newcommand{\FF}{\mathbb{F}}
\newcommand{\RR}{\mathbb{R}}
\newcommand{\ZZ}{\mathbb{Z}}
\newcommand{\PP}{\mathbb{P}}
\newcommand{\TT}{\mathbb{T}}
\newcommand{\HH}{\mathbb{H}}
\newcommand{\NN}{\mathbb{N}}

\newcommand{\Hbar}{{\overline{\mathbb{H}}}}

\newcommand{\Fbar}{{\overline{\FF}}}
\newcommand{\res}{\mathrm{res}}
\newcommand{\cts}{\mathrm{cts}}
\newcommand{\Mod}{{\mathfrak{Mod}}}
\newcommand{\Shf}{{\mathfrak{Sh}}}
\newcommand{\Hom}{{\rm Hom}}
\newcommand{\SL}{\mathrm{SL}}
\newcommand{\cusps}{\mathrm{cusps}}
\DeclareMathOperator{\Image}{im}
\DeclareMathOperator{\Ker}{ker}
\newcommand{\cE}{\mathcal{E}}
\newcommand{\cC}{\mathcal{C}}
\newcommand{\cB}{\mathcal{B}}
\newcommand{\cI}{\mathcal{I}}
\newcommand{\cO}{\mathcal{O}}
\newcommand{\cF}{\mathcal{F}}
\newcommand{\cV}{\mathcal{V}}
\newcommand{\calH}{\mathcal{H}}
\newcommand{\cM}{\mathcal{M}}
\newcommand{\cCM}{\mathcal{CM}}
\newcommand{\Ind}{{\rm Ind}}
\newcommand{\Res}{{\rm Res}}
\newcommand{\Coind}{{\rm Coind}}
\newcommand{\PSL}{\mathrm{PSL}}
\newcommand{\mat}[4]{
 \left(  \begin{smallmatrix} #1 & #2 \\ #3 & #4 \end{smallmatrix} \right)}
\newcommand{\Hpar}{H_{\mathrm{par}}}
\newcommand{\Hc}{H_{\mathrm{c}}}
\newcommand{\ilim}{{\underset{\to}{\lim}}}
\newcommand{\pf}{{\bf Proof. }}
\newcommand{\qed}{\hspace* {.5cm} \hfill $\Box$}

\begin{document}

\title{On modular symbols and the cohomology of\\ Hecke triangle surfaces}
\author{Gabor Wiese}
\maketitle

\begin{abstract}
The aim of this article is to give a concise algebraic treatment
of the modular symbols formalism, generalized
from modular curves to Hecke triangle surfaces.
A sketch is included of how the modular symbols formalism gives rise
to the standard algorithms for the computation of holomorphic modular forms.
Precise and explicit connections are established to the 
cohomology of Hecke triangle surfaces and group cohomology.
In all the note a general commutative ring is used as coefficient
ring in view of applications to the computation of modular forms
over rings different from the complex numbers.

\end{abstract}

\section{Introduction}

The purpose of this article is to give a concise algebraic treatment
of the modular symbols formalism and related objects 
over an arbitrary commutative ring.
We show how the standard algorithms for computing modular forms
for congruence subgroups of $\SL_2(\ZZ)$ can be deduced purely
algebraically, using only the Eichler-Shimura isomorphism.
This includes an algebraic proof of the presentation of modular symbols
in terms of Manin symbols.
Thus, we avoid the use of the technically rather difficult paper
\cite{Shokurov}, which is present in all the published treatments
known to the author (e.g.\ \cite{MerelUniversal}, \cite{SteinBook}).

The algebraic formulation that we give generalizes immediately to the
so-called Hecke triangle groups and we choose that general
set-up from the beginning.
It should, however, be pointed out that most Hecke triangle groups
are non-arithmetic. In that case our treatment only gives an
isomorphism between generalized modular symbols and modular forms,
but does not yield an algorithm for the computation of the Fourier
coefficients of the modular forms due to the absence of a suitable
Hecke theory.

There is considerable interest in trying to use the modular symbols
formalism over rings different from the complex number in order to
compute modular forms over these rings (e.g.\ using the theory
of Katz modular forms, see \cite{EdixJussieu} or \cite{Thesis}).
This is the reason why the treatment of the modular symbols formalism
in the present article is over any commutative unitary ring. In this
generality (in fact already over the integers) one notices quickly
that, although the modular symbols formalism is inspired by the homology
of modular curves, it does not quite compute it (compare 
Theorem~\ref{ManinSymbols} and Remark~\ref{RemMerel}).

Moreover, in order to be able to treat questions like the computation
of modular forms over more general rings one needs a geometric or
algebraic interpretation of the modular symbols formalism
in order to be able to establish a link with modular forms.
For this reason we also treat certain cohomology groups on
modular curves (as Riemann surfaces) and certain group cohomology
groups which are both closely related to the modular symbols formalism.
All the objects appearing are described by explicit formulae and
a precise comparison is included. 
The cohomology group of the modular curve considered in the present
article is a complex analogue of the \'etale cohomology group
that one uses to define the $2$-dimensional $l$-adic Galois representations
attached to a modular form. 

The differences between the various objects come from non-trivially
stabilized points on the upper half plane.
It is hence natural to use analytic modular
stacks instead of modular curves and compare these two
via the Leray spectral sequence.
This has been carried out in the author's thesis. However, for the sake
of the present article a formulation was chosen that uses only
homological algebra, but gives the same results.

Apart from some facts about Hecke groups and some cohomology theory of
groups and topological spaces, the treatment of the article is
essentially self-contained.

\subsection*{Overview}

In Section~\ref{SecHtg} we present some facts about Hecke triangle
groups, as well as two results to be used in the sequel.
The following three sections are independent of one another.
In Section~\ref{MSSec} we introduce
the modular symbols formalism extended to subgroups of finite index in
Hecke triangle groups, and give a description in terms of Manin
symbols.  
An explicit formula for the parabolic subspace of the group cohomology
for subgroups of finite index of Hecke triangle groups is
derived in Section~\ref{GroupSec}. 
The subsequent Section~\ref{CurveSec} treats a similar
cohomology group for a certain sheaf on the modular surface
for~$\Gamma$. An explicit formula is derived, which generalizes a
result of Merel's to higher weights. In the final section a comparison
between the objects is carried out and it is sketched how the
Eichler-Shimura theorem together with a theory of Hecke operators and
the results from the previous sections can be used to compute modular
forms for congruence subgroups of $\SL_2(\ZZ)$.

\subsection*{Acknowledgements}

This article has grown out of Chapter~II of my thesis.  I would like
to thank my thesis advisor Bas Edixhoven for many useful suggestions,
and my office mate Theo van den Bogaart for sharing his insights in
algebraic geometry, and William Stein for comments on an early
version.

\section{Hecke triangle groups and surfaces}\label{SecHtg}

For an integer $n \ge 3$ one defines the $n$-th 
{\em Hecke triangle group}~$\Delta_n$ as the subgroup
of $\PSL_2(\RR)$ generated by
$$ \sigma := \mat 0 {-1} 1 0 \;\;\; \text{ and } \;\;\;
  \tau := \mat {\lambda_n} {-1} 1 0$$
with $\lambda_n = e^{\pi i/n} + e^{-\pi i/n} = 2 \cos(\pi/n)$.
In abuse of notation when writing a matrix we often mean its class
modulo scalar matrices, i.e.\ an element of the projective
linear group.
The generation is free and $\Delta_n$ is the free product
of $\ZZ/2\ZZ = \langle \sigma \rangle$ and $\ZZ/n\ZZ = \langle \tau \rangle$, 
which makes the cohomological computations to come very simple.
The stabilizer of the element $\infty \in \PP^1(\RR)$ is 
$(\Delta_n)_\infty = \langle T \rangle$ with $T = \tau \sigma = \mat 1 {\lambda_n} 0 1$.
We denote by $\Delta_n(\infty)$ the orbit of $\infty$ under~$\Delta_n$.
As a special case let us note that $\Delta_3 = \PSL_2(\ZZ)$ and $\Delta_3(\infty) = \PP^1(\QQ)$.
The Hecke group $\Delta_n$ is a Fuchsian group of the first kind
\begin{wrapfigure}{r}{6cm}\label{fig}
\includegraphics[width=6cm]{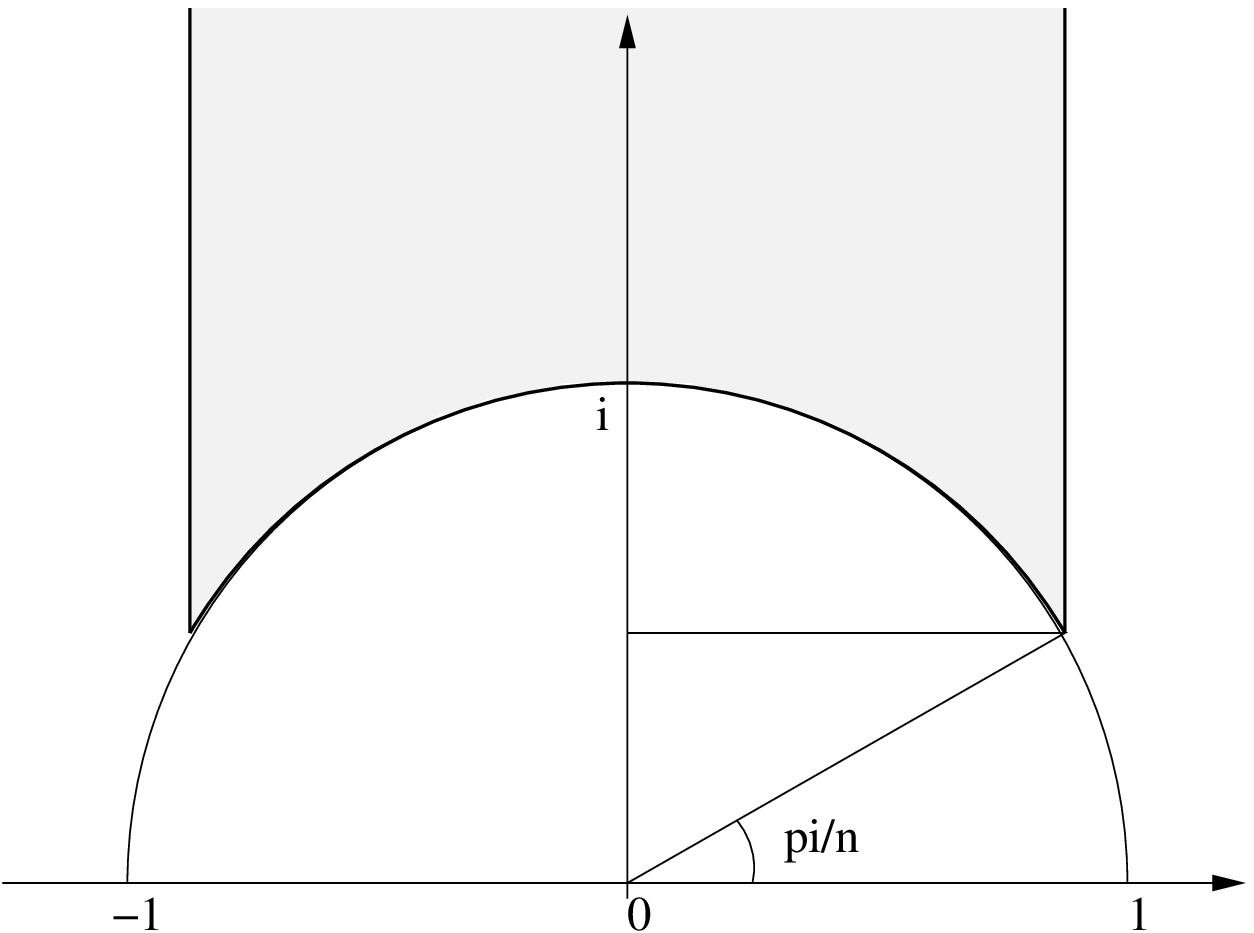}
\end{wrapfigure}
(i.e.\ a discrete subgroup of~$\PSL_2(\RR)$ of finite
covolume) having parabolic elements (namely precisely the conjugates of powers of~$T$).
Any subgroup $\Gamma \le \Delta_n$ of finite index is also a Fuchsian group
of the first kind with parabolic elements.
Moreover, $Y_\Gamma := \Gamma \backslash \HH$ can be given
the structure of an open Riemann surface.
It can be compactified $X_\Gamma := Y_\Gamma \cup \Gamma \backslash \Delta_n(\infty)$,
where the set $\Gamma \backslash \Delta_n(\infty)$ is the set of
{\em (parabolic) cusps of~$\Gamma$}. We write $\Hbar = \HH \cup \Delta_n(\infty)$.
The compact Riemann surface $X_{\Delta_n}$ is called a {\em Hecke triangle surface}.
The open Riemann surface $Y_{\Delta_n}$ can be visualized as a fundamental domain with angles
$0$, $\pi/n$ and~$\pi/n$, as shown in Fig.~\ref{fig}.

By a result of Leutbecher (as cited in \cite{Schmidt})
Hecke triangle groups are non-arithmetic, except for $n \in \{3,4,6\}$,
i.e.\ $\Delta_n$ is not commensurable with any $\PSL_2(\cO)$ for
$\cO$ the ring of integers of a number field. For more details on Hecke triangle
groups and surfaces we refer to \cite{Schmidt} and the references therein.

The choice of working inside projective linear groups instead of linear groups
was made since it simplifies some formulae and some proofs at nearly no costs.

\subsection*{Notation}

In most of the paper we use the following notation.

\begin{notation}\label{notation}
Let $G = \Delta_n < \PSL_2(\RR)$ be a Hecke triangle group for some integer $n \ge 3$
and let $\Gamma \le G$ be a subgroup of finite index.
We use the notations $Y_\Gamma$, $X_\Gamma$ for the Riemann surfaces
introduced above and call $j: Y_\Gamma \hookrightarrow X_\Gamma$
the natural embedding and
$\pi: \HH   \twoheadrightarrow Y_\Gamma$ resp.\ 
$\pi: \Hbar \twoheadrightarrow X_\Gamma$
the natural projections. 
Furthermore, the matrices $\sigma, \tau, T$ defined above will be used.

We let $R$ be a commutative ring with unit and $V$ a left $R[\Gamma]$-module.
If $g \in G$ is some element of finite order~$m$, we denote by $N_g$ the element
$1 + g + \dots + g^{m-1}$ of the group ring $R[G]$. Similarly, if $H \le G$
is a finite subgroup, we write $N_H = \sum_{h \in H} h \in R[G]$.
\end{notation}

\subsection*{Mayer-Vietoris and amalgamated products}

We assume Notation~\ref{notation}.

\begin{proposition}\label{mayervietoris}
Let $M$ be a left $R[G]$-module. Then the Mayer-Vietoris
sequence gives the exact sequences
\begin{align*}
0 \to M^{G} \to M^{\langle \sigma \rangle} \oplus M^{\langle \tau \rangle} \to M \to 
H^1(G,M) \to H^1(\langle \sigma \rangle, M) \oplus H^1(\langle \tau \rangle, M) \to 0,\\
0 \to H_1(\langle \sigma \rangle, M) \oplus H_1(\langle \tau \rangle, M) \to
H_1(G,M) \to M \to M_{\langle \sigma \rangle} \oplus M_{\langle \tau \rangle} \to M_G \to 0
\end{align*}
and for all $i \ge 2$ isomorphisms
\begin{align*}
H^i(G,M) \cong H^i(\langle \sigma \rangle, M) \oplus H^i(\langle \tau \rangle, M),\\
H_i(G,M) \cong H_i(\langle \sigma \rangle, M) \oplus H_i(\langle \tau \rangle, M).\\
\end{align*}
\end{proposition}

\pf
Let us write $G_1 := \langle \sigma \rangle$ and $G_2 := \langle \tau \rangle$. 
By \cite{Brown}, II.8.8, we have the exact sequence
$$ 0 \to R[G] \to R[G/G_1] \oplus R[G/G_2] \to R \to 0$$
of $R[G]$-modules, which are free as $R$-modules.
Application of the functor 
$\Hom_R(\cdot, M)$ gives rise to the exact sequence of $R[G]$-modules
$$ 0 \to M \to \Hom_{R[G_1]}(R[G],M) \oplus \Hom_{R[G_2]}(R[G],M) \to
\Hom_R (R[G],M) \to 0.$$
The central terms, as well as the term on the right, can be identified
with coinduced modules.
Hence, the statements on cohomology follow by taking the long exact sequence of cohomology
and invoking Shapiro's lemma.
Using the functor $\cdot \otimes_R M$ gives rise to the analogous statements
about homology.
\qed

\subsection*{Mackey's formula and stabilizers}

We now prove Mackey's formula for coinduced modules.
If $H \le G$ are groups and $V$ is an $R[H]$-module, 
the coinduced module $\Coind_H^{G} V$ 
can be described as $\Hom_{R[H]}(R[G],V)$.

\begin{proposition}\label{propmackey}
Let $R$ be a ring, $G$ a group and $H,K$ subgroups of~$G$. 
Let furthermore $V$ be an $R[H]$-module. Then {\em Mackey's formula}
$$ \Res_K^G \Coind_H^G V \cong \prod_{g \in H\backslash G / K}
\Coind_{K\cap g^{-1}Hg}^K {}^g (\Res^H_{H\cap gKg^{-1}} V)$$
holds. Here ${}^g (\Res^H_{H\cap gKg^{-1}} V)$ denotes the
$R[K \cap g^{-1}Hg]$-module obtained from $V$ via the 
conjugated action $g^{-1}hg ._g v := h. v$ for $v \in V$
and $h \in H$ such that $g^{-1}hg \in K$.
\end{proposition}

\pf
We consider the commutative diagram
$$ \xymatrix@=0.9cm{
\Res_K^G \Hom_H(R[G],V) \ar@{->}[r]  \ar@{->}[dr] &
\prod_{g \in H\backslash G / K}
\Hom_{K\cap g^{-1}Hg}(R[K], {}^g (\Res^H_{H\cap gKg^{-1}} V)) \ar@{->}^\sim[d] \\
&\prod_{g \in H\backslash G / K}
 \Hom_{H\cap gKg^{-1}}(R[gKg^{-1}], \Res^H_{H\cap gKg^{-1}} V)).}$$
The vertical arrow is just given by conjugation and is clearly
an isomorphism.
The diagonal map is the product of the natural restrictions.
From the bijection
$$ \big(H \cap gKg^{-1}\big) \backslash gKg^{-1} 
\xrightarrow{gkg^{-1} \mapsto Hgk} H \backslash HgK$$
it is clear that also the diagonal map is an isomorphism,
proving the proposition.
\qed

Applying Shapiro's lemma, one immediately obtains the following 
two corollaries.

\begin{corollary}\label{cormackey}
In the situation of Proposition~\ref{propmackey} one has
\begin{align*}
H^i(K,\Coind_H^G V) 
& \cong \prod_{g \in H\backslash G / K}
H^i(K \cap g^{-1}Hg,  {}^g (\Res^H_{H\cap gKg^{-1}} V) \\
& \cong \prod_{g \in H\backslash G / K}
H^i(H \cap gKg^{-1},  \Res^H_{H \cap gKg^{-1}} V) 
\end{align*}
for all $i \in \NN$.
\end{corollary}

\begin{corollary}\label{corstabmackey}
We now assume Notation~\ref{notation}. For $x \in \Hbar$ we denote
by $G_x$ the stabilizer subgroup of~$G$ for the point~$x$. 
The image of the $G$-orbit of $x$ in $X_\Gamma$ is in bijection with
the double cosets $\Gamma \backslash G / G_x$ as follows 
$$ \Gamma \backslash G / G_x \xrightarrow{g \mapsto gx} \Gamma \backslash Gx.$$
Moreover, the group $\Gamma \cap gG_x g^{-1}$ equals
$\Gamma_{gx}$, the stabilizer subgroup of~$\Gamma$ for the point~$gx$.
Thus, for all $i \in \NN$, Mackey's formula gives an isomorphism
$$H^i(G_x, \Coind_\Gamma^G V) \cong \prod_{y \in \Gamma \backslash Gx} 
H^i(\Gamma_y, V).$$
\end{corollary}

\section{The modular symbols formalism}\label{MSSec}

Modular symbols were systematically studied by Manin(\cite{Manin}).
In Cremona's book \cite{Cremona} 
it is shown how the modular symbols formalism can be used
for computing weight two modular forms. A generalization to higher weight
modular forms was found by \u{S}okurov (\cite{Shokurov}),
and the resulting ``higher weights modular symbols formalism'' was
used by Merel (\cite{MerelUniversal}) algorithmically. 
This formalism is also one of the subjects of William Stein's very 
comprehensive textbook \cite{SteinBook}.

Common to the mentioned published treatments is that they use rather difficult
computations of homology groups.  
We show in this section that the description of modular symbols in
terms of Manin symbols can be derived purely algebraically with simple
and conceptual methods.  Together with the Eichler-Shimura
isomorphism, the proof of which is also quite easy, this already gives the
basic idea of the algorithms of \cite{SteinBook}. This is sketched in
Section~\ref{MFSec}. These algorithms were implemented in {\sc Magma} and
{\sc Sage} by William Stein in the case of the standard congruence
subgroups.

\subsection*{Definition}

Modular symbols can be thought of as geodesic paths between two cusps
resp.\ as the associated homology class relative to the cusps. We shall,
however, give a combinatorial definition, as is implemented in {\sc Magma}
and like the one in \cite{MerelUniversal}, \cite{Cremona} 
and~\cite{SteinBook}, 
except that we do not factor out torsion, 
but intend a common treatment for all rings. 

We give the definition in the more general context of Hecke
triangle groups and also allow general modules.
Throughout this section we assume Notation~\ref{notation}.

\begin{definition}\label{defMS}
We define the $R$-modules
$$ \cM_R := R[\{\alpha,\beta\}| \alpha,\beta \in G(\infty)]/
\langle \{\alpha,\alpha\}, \{\alpha,\beta\} + \{\beta,\gamma\} + \{\gamma,\alpha\}
| \alpha,\beta,\gamma \in G(\infty)\rangle$$
and
$$ \cB_R := R[G(\infty)].$$ 
We equip both with the natural left $\Gamma$-action. 
Furthermore, we let
$$ \cM_R(V) := \cM_R \otimes_R V \;\;\;\; \text{ and } \;\;\;\; \cB_R(V) := \cB_R \otimes_R V$$
with the left diagonal $\Gamma$-action.

\begin{enumerate}
\item We call the $\Gamma$-coinvariants
$$ \cM_R (\Gamma,V) :=  \cM_R(V)_\Gamma = 
\cM_R(V)/ \langle (x - g x) | g \in \Gamma, x \in \cM_R(V) \rangle$$ 
{\em the space of $(\Gamma,V)$-modular symbols.}

\item We call the $\Gamma$-coinvariants
$$ \cB_R(\Gamma,V) :=  \cB_R(V)_\Gamma = 
\cB_R(V)/ \langle (x - g x) | g \in \Gamma, x \in \cB_R(V) \rangle$$ 
{\em the space of $(\Gamma,V)$-boundary symbols.}

\item We define the {\em boundary map} as the map
$$ \cM_R(\Gamma,V) \to \cB_R(\Gamma,V)$$
which is induced from the map $\cM_R \to \cB_R$ sending $\{\alpha, \beta\}$
to $\{\beta\} - \{\alpha\}$.

\item The kernel of the boundary map is denoted by $\cCM_R(\Gamma,V)$ and is called 
{\em the space of cuspidal $(\Gamma,V)$-modular symbols.}

\item The image of the boundary map inside $\cB_R(\Gamma,V)$ is
denoted by $\cE_R(\Gamma,V)$ and is called
{\em the space of $(\Gamma,V)$-Eisenstein symbols.}
\end{enumerate}
\end{definition}

\subsection*{Manin symbols}

Manin symbols provide an explicit description of modular symbols. 
We stay in the general setting over a ring~$R$ and keep 
Notation~\ref{notation}.

As $G$ is infinite, the induced module $R[G]$ is not isomorphic to
the coinduced one $\Hom_R(R[G],R)$ and $R[G]$ is not cohomologically
trivial. However, $\Hpar^1(G,R[G])=0$ (for the definition see
Section~\ref{GroupSec}). This is the essence of the following proposition.

\begin{proposition}\label{hparnull}
The sequence of $R$-modules
$$ 0 \to R[G]N_\sigma + R[G]N_\tau \to R[G]
\xrightarrow{g \,\mapsto\, g(1-\sigma)\infty} R[G(\infty)] 
\xrightarrow{g\infty \,\mapsto\, 1} R \to 0$$
is exact.
\end{proposition}

\pf
We first use that $R[G]$ is a cohomologically trivial module for 
both $\langle \sigma \rangle$ and $\langle \tau \rangle$. This gives
$$ R[G]N_\sigma = \ker_{R[G]} (1-\sigma) = R[G]^{\langle \sigma \rangle}, \;\;\; 
R[G]N_\tau = \ker_{R[G]} (1-\tau)= R[G]^{\langle \tau \rangle}, $$
$$ R[G](1-\sigma) = \ker_{R[G]} N_\sigma \;\;\; \text{ and } 
\;\;\; R[G](1-\tau) = \ker_{R[G]} N_\tau.$$
Proposition~\ref{mayervietoris} yields the exact sequence
$$ 0 \to R[G] \to R[G]_{\langle \sigma \rangle} \oplus R[G]_{\langle \tau \rangle}
 \to R \to 0,$$
since $H_1(G,R[G]) = 0$ by Shapiro's lemma because $R[G] \cong \Ind_1^G R$.
In fact, this sequence is at the origin of our proof of 
Proposition~\ref{mayervietoris}.
The injectivity of the first map in the exact sequence means
\begin{equation}\label{kereq}
R[G](1-\sigma) \cap R[G](1-\tau) = 0.
\end{equation}

We identify $R[G]/R[G](1-T)$ with $R[G(\infty)]$ by
sending $g$ to $g\infty$.
Now we show the exactness at $R[G]$,
which comes down to proving that the equation
$ x(1-\sigma) = y (1-T)$ for $x,y \in R[G]$ implies that
$x$ is in $R[G]^{\langle \sigma \rangle} + R[G]^{\langle \tau \rangle}$.

Using the formula $\tau = T\sigma$ we obtain that
$x(1-\sigma) = y(1-T) = y(1-\tau) -yT(1-\sigma)$.
This yields $x(1-\sigma) +yT(1-\sigma) = y(1-\tau)$. This expression, however,
is zero by Eq.~(\ref{kereq}). 
Consequently, there is a $z \in R[G]$ such that $y = zN_\tau$.
Hence, using $T = \tau \sigma$ and consequently $N_\tau T = N_\tau \sigma$, we get
$$ y(1-T) = z N_\tau(1-T) = z N_\tau (1-\sigma) = y(1-\sigma).$$
The equation $x(1-\sigma) = y(1-\sigma)$ means that 
$x-y$ is in $R[G]^{\langle \sigma \rangle}$. 
As we know that $y \in R[G]^{\langle \tau \rangle}$,
we see that $x = (x-y) + y$ is in $R[G]^{\langle \sigma \rangle} + R[G]^{\langle \tau \rangle}$,
as required.
Note that instead of this explicit calculation we could also have appealed
to Proposition~\ref{propPara}.

The exactness at $R[G(\infty)]$ can be seen as follows (we avoid here
the traditional continued fractions argument, as it does not obviously
generalize to Hecke triangle groups and is not in the spirit of the
present group theoretic approach).
Since $\sigma$ and $T=\tau\sigma$ generate~$G$, the kernel of
$R[G] \xrightarrow{g \mapsto 1} R$ is
$R[G](1-\sigma) + R[G](1-T)$. Taking the quotient by $R[G](1-T)$
gives the desired exactness.
\qed

\begin{lemma}\label{mrlem}
The sequence of $R$-modules
$$ 0 \to \cM_R \xrightarrow{\{\alpha,\beta\} \mapsto \beta - \alpha} R[G(\infty)] 
\xrightarrow{\alpha \mapsto 1} R \to 0$$
is exact.
\end{lemma}

\pf
The injectivity of the first arrow is clear, since 
we can write any element in
$\cM_R$ as $\sum_{\alpha \neq \infty} r_\alpha \{\infty, \alpha\}$ with $r_\alpha \in R$,
using the relations defining~$\cM_R$. The image of this element under the first arrow
is $\sum_{\alpha \neq \infty} r_\alpha \alpha - (\sum_{\alpha \neq \infty} r_\alpha) \infty$.
If this is zero, clearly all $r_\alpha$ are zero, proving the injectivity of the first
arrow.

Suppose now we are given $\sum_\alpha r_\alpha \alpha \in R[G(\infty)]$ in the kernel
of the second arrow. Then $\sum_\alpha r_\alpha = 0$ and consequently we have
$$\sum_\alpha r_\alpha \alpha = \sum_{\alpha \neq \infty} r_\alpha \alpha - 
(\sum_{\alpha \neq \infty} r_\alpha) \infty$$
which is in the image of the first arrow, as noticed before.
\qed

\begin{proposition}\label{propker}
The homomorphism of $R$-modules
$$ R[G] \xrightarrow{\phi} \cM_R,\;\;\;
g \mapsto \{g.0,g.\infty\}$$
is surjective and its kernel is given by 
$R[G]N_\sigma + R[G]N_\tau$. 
\end{proposition}

\pf
This is a direct consequence of Proposition~\ref{hparnull}
and Lemma~\ref{mrlem}.
\qed

We are now ready to prove the description of modular symbols in terms of
Manin symbols.

\begin{theorem}\label{ManinSymbols}
Recall that we are assuming Notation~\ref{notation}.
Let $M = \Ind_\Gamma^{G} V$, which we identify with
$(R[G] \otimes_R V)_\Gamma$. 
That module carries the right $R[G]$-action
$(h \otimes v)g = (hg \otimes v)$ for $g,h \in G$, $v \in V$,
and the $\Gamma$-coinvariants are taken
for the diagonal left $\Gamma$-action.
The following statements hold:
\begin{enumerate}
\item\label{ma}
The homomorphism $\phi$ from Proposition~\ref{propker} induces the
exact sequence of $R$-modules
$$ 0 \to M N_\sigma + M N_\tau \to M \to \cM_R(\Gamma,V) \to 0.$$
\item\label{mb}
The homomorphism $R[G] \to R[G(\infty)]$ sending $g$ to $g.\infty$
induces the exact sequence of $R$-modules
$$ 0 \to M(1-T) \to M \to \cB_R(\Gamma,V) \to 0.$$
\item\label{mc}
The identifications of (\ref{ma}) and~(\ref{mb}) imply the isomorphism
$$ \cCM_R(\Gamma,V) \cong \ker\big(
M/ (M N_\sigma + M N_\tau ) \xrightarrow{m \mapsto m(1-\sigma)} M/M(1-T) \big).$$
\end{enumerate}
\end{theorem}

\pf
(\ref{ma}) We derive this from Proposition~\ref{propker}, which gives
the exact sequence
$$ 0 \to R[G]N_\sigma + R[G]N_\tau \to R[G] \to \cM_R \to 0.$$
Tensoring with $V$ over $R$, we obtain the exact sequence of 
left $R[\Gamma]$-modules
$$ 0 \to (R[G] \otimes_R V) N_\sigma +  (R[G] \otimes_R V) N_\tau \to (R[G] \otimes_R V) \to \cM_R(V) \to 0.$$
Passing to left $\Gamma$-coinvariants yields~(\ref{ma}). 
Part~(\ref{mb}) is clear from the definition and Part~(\ref{mc}) has already been noticed
in the proof of Proposition~\ref{hparnull}.
\qed

In the literature on Manin symbols one usually finds different
versions of the module~$M$, namely the following.
Suppose first that the $\Gamma$-action on~$V$ is the restriction of some $G$-action on~$V$.
Then we have the isomorphism
$$  (R[G] \otimes_R V)_\Gamma \cong R[\Gamma \backslash G] \otimes_R V, \;\;\; 
g \otimes v \mapsto g \otimes g^{-1} v.$$
The right $R[G]$-action carries over to the action
$(\Gamma h \otimes v)g = \Gamma hg \otimes g^{-1} v$.

We should also mention a slight variant of this. Suppose now that $\tilde{G}$
and $\tilde{\Gamma}$ are defined as $\phi^{-1}(G)$ resp.\ $\phi^{-1}(\Gamma)$
for the projection $\phi: \SL_2(\RR) \twoheadrightarrow \PSL_2(\RR)$.
We also assume that $-1 \in \tilde{\Gamma}$, so that $V$ is an $R[\tilde{\Gamma}]$-module,
and that the $\tilde{\Gamma}$-action on~$V$ is the restriction of some $\tilde{G}$-action on~$V$.
Then we have the isomorphism
$$ (R[G] \otimes_R V)_\Gamma \cong R[\tilde{\Gamma} \backslash \tilde{G}] \otimes_R V, \;\;\; 
g \otimes v \mapsto g \otimes g^{-1} v.$$

\section{Group cohomology}\label{GroupSec}

Also in this section we assume Notation~\ref{notation}.

\subsection*{Definitions}

We define {\em parabolic group cohomology} as the
left hand term and the {\em boundary  group cohomology}
as the right hand term in the exact sequence
$$ 0 \to \Hpar^1 (\Gamma, V) \to H^1(\Gamma, V) \xrightarrow{\res} 
\prod_{c \in \Gamma \backslash G(\infty)} 
H^1(\Gamma_{\tilde{c}}, V),$$
where $\Gamma_{\tilde{c}}$ is the stabilizer subgroup of~$\Gamma$
for the cusp~$\tilde{c} \in \Hbar$ with $\pi(\tilde{c}) = c$.
We point out that (parabolic) group cohomology would in general
be different if we worked with subgroups of $\SL_2(\RR)$
and not $\PSL_2(\RR)$ throughout. 

\subsection*{Computing group cohomology}

In order to compute the group cohomology for~$\Gamma$,
it suffices to compute the cohomology for~$G$
because of Shapiro's lemma,
which for any $R[\Gamma]$-module~$V$ gives an isomorphism
$$H^1(G, \Coind_\Gamma^{G} V) \cong H^1(\Gamma, V).$$
Due to Corollary~\ref{corstabmackey} it is clear
that Shapiro's lemma respects the parabolic subspace.

A first, however, not complete computation of the
group cohomology of $R[G]$-modules is provided by
the Mayer-Vietoris sequence (Proposition~\ref{mayervietoris}).
We now derive an explicit description.

\begin{proposition}\label{propcoh}
Let $M$ be a left $R[G]$-module. Then the sequence of $R$-modules
$$ 0 \to M^{G} \to M \to \ker_M N_\sigma \times \ker_M N_\tau 
\to H^1(G,M) \to 0$$
is exact.
\end{proposition}

\pf
We determine the $1$-cocycles of~$M$. Apart from $f(1) = 0$, they must satisfy
$$ 0 = f(\sigma^2) = \sigma f(\sigma) + f(\sigma) = N_\sigma f(\sigma) \text{ and}$$
$$ 0 = f(\tau^n) = \dots = N_\tau f(\tau).$$
Since these are the only relations in $G$, 
a cocycle is uniquely given by the choices
$$ f(\sigma) \in \ker_M N_\sigma \text{ and }
   f(\tau) \in \ker_M N_\tau.$$
The $1$-coboundaries are precisely the cocycles $f$ which satisfy 
$f(\sigma) = (1-\sigma)m$ and $f(\tau)=(1-\tau)m$ for some $m \in M$.
This proves
$$ H^1(G,M) \cong (\ker_M N_\sigma \times \ker_M N_\tau) /
\big( ((1-\sigma)m,(1-\tau)m) \; | \, m \in M \big).$$
Rewriting yields the proposition. 
\qed

\begin{remark}\label{remU}
As $G_\infty = \langle T \rangle < G$ is infinite cyclic, 
one has $H^1(G_\infty,\Res_{G_\infty}^G M) \cong M / (1-T)M$.
\end{remark}

An explicit presentation of the parabolic group cohomology is the following. 

\begin{proposition}\label{propPara}
The parabolic group cohomology group sits in the exact sequence
$$ 0 \to M^{\langle T \rangle}/M^{G} \to \ker_M N_\sigma \cap \ker_M N_\tau 
\xrightarrow{\phi} \Hpar^1(G,M) \to 0,$$
where $\phi$ maps an element~$m$ to the $1$-cocycle~$f$ uniquely determined
by $f(\sigma) = f(\tau) = m$. 
\end{proposition}

\pf
Using Proposition~\ref{propcoh} we get the exact commutative diagram
$$ \xymatrix@=1.3cm{
M^{\langle T \rangle}/M^{G} \ar@{^{(}->}^(.47){(\sigma^{-1} - 1)}[r] \ar@{^{(}->}^{\sigma^{-1}}[d] & 
\ker N_\sigma \cap \ker N_\tau \ar@{->}[r] \ar@{^{(}->}[d] &
\Hpar^1(G,M) \ar@{^{(}->}[d] \\
M/M^{G} \ar@{^{(}->}^(.46){(1-\sigma, 1-\tau)}[r] \ar@{->>}^{(1-T)\sigma}[d] & 
\ker N_\sigma \times \ker N_\tau \ar@{->>}[r] \ar@{->}^{(a,b)\mapsto b-a}[d] &
H^1(G,M) \ar@{->}[d] \\
(1-T)M \ar@{^{(}->}[r] & 
M \ar@{->>}[r] &
H^1(G_\infty,M).  }$$
As the bottom left vertical arrow is surjective, the claim follows from
the snake lemma.
\qed

\section{Cohomology of Hecke triangle surfaces}\label{CurveSec}

\subsection*{The group cohomology presheaf and sheaf}

In this section we let $X$ be a topological space, $R$ a commutative ring
and $\Gamma$ a group.
For any ring $S$, not necessarily commutative, we denote by
$S-\Mod$ the category of left $S$-modules and by
$\Shf_X(S-\Mod)$ the category of sheaves of left $S$-modules on~$X$.

We collect some well-known, but important, properties
in the following proposition.

\begin{proposition}\label{propss}
\begin{enumerate}
\item \label{pa}
The category $\Shf_X(S-\Mod)$ has enough injectives.
\item \label{pb}
Let $\cI \in \Shf_X(S-\Mod)$ be an injective object.
Then $\cI$ is flabby ({\it flasque}).
\item \label{pc}
Let $\cV$ be an object of $\Shf_X(R[\Gamma]-\Mod)$.
Then the cohomology groups $H^i(X,\cV)$ for $i \ge 0$ do not depend
on whether they are computed in the category $\Shf_X(R[\Gamma]-\Mod)$
or by forgetting the $\Gamma$-action in the category $\Shf_X(R-\Mod)$.
\item \label{pd}
Let $\cI \in \Shf_X(S-\Mod)$ be an injective object. Then for
all open sets $U \subseteq X$, the object $\cI(U)$ of $S-\Mod$ is
injective.
\item \label{pe}
Let $\cI \in \Shf_X(R[\Gamma]-\Mod)$ be an injective object.
Then $\cI^\Gamma$ is an injective object of $\Shf_X(R-\Mod)$.
\end{enumerate}
\end{proposition}

\pf
We notice that $X$ together with the constant sheaf $S$ on $X$
is a ringed space. The statements (\ref{pa}) and (\ref{pb}) are then
\cite{Hartshorne}, Proposition~III.2.2 and Lemma~III.2.4.
We should, however, point out that Hartshorne works with
commutative rings only. But the proofs also work in the
non-commutative situation.

(\ref{pc}) follows from (\ref{pb}), as flabby resolutions can be used for
computing $H^i(X,\cV)$ (see \cite{Hartshorne}, Proposition~III.2.5).

(\ref{pd}) As the sheaf $\cI$ restricted to $U$ is injective, it
suffices to prove the statement for $U=X$.
The injectivity of $\cI$ means that the functor
$\Hom_{\Shf_X(S-\Mod)}(\cdot,\cI)$ is exact.
For $A \in S-\Mod$, we denote by $\underline{A}$ the
constant sheaf on~$X$ associated with~$A$. We have
$$ \Hom_{S-\Mod}(A, \cI(X)) = \Hom_{\Shf_X(S-\Mod)}(\underline{A},\cI)$$
for all $A \in S-\Mod$. As taking the constant sheaf is an exact
functor, $\Hom_{S-\Mod}(\cdot, \cI(X))$ is also exact, proving
the injectivity of $\cI(X)$.

(\ref{pe}) Suppose we are given a diagram 
$$ \xymatrix@=.5cm{
& \cI^\Gamma \\
A  \ar@{^{(}->}[r] \ar@{->}[ur] & B }$$
in $\Shf_X(R-\Mod)$. By composing with the natural injection
$\cI^\Gamma \hookrightarrow \cI$ and putting a trivial $\Gamma$-action
on $A$ and $B$, we obtain the commutative diagram
$$ \xymatrix@=.5cm{
& \cI \\
A  \ar@{^{(}->}[r] \ar@{->}[ur] & B \ar@{->}[u]}$$
in $\Shf_X(R[\Gamma]-\Mod)$, since $\cI$ is injective.
However, the image of $B \to \cI$ is contained in $\cI^\Gamma$,
as $B$ is a trivial $\Gamma$-module.
\qed

Let $U \subseteq X$ be an open set. We consider the following
commutative diagram of categories:
$$ \xymatrix@=.7cm{
& \Shf_X(R-\Mod) \ar@{->}[dr]^{H^0(U,\cdot)} & \\
\Shf_X(R[\Gamma]-\Mod) \ar@{->}[ur]^{(\cdot)^\Gamma} \ar@{->}[dr]^{H^0(U,\cdot)} && R-\Mod.\\
& R[\Gamma]-\Mod \ar@{->}[ur]^{(\cdot)^\Gamma}& }$$
Let $\cV \in \Shf_X(R[\Gamma]-\Mod)$.
Due to Proposition~\ref{propss},
Grothendieck's theorem on spectral sequences (\cite{Weibel}, Theorem~5.8.3)
gives rise to the two spectral sequences
\begin{equation}\label{spseins}
H^p(\Gamma,H^q(U,\cV)) \Rightarrow R^{p+q}(H^0(U,\cdot)^\Gamma)(\cV)
\end{equation}
and
\begin{equation}\label{spszwei}
H^p(U,\calH^q(\Gamma,\cV))\Rightarrow R^{p+q}(H^0(U,\cdot)^\Gamma)(\cV), 
\end{equation}
where we write 
$$\calH^q(\Gamma,\cV) = R^q((\cdot)^\Gamma)(\cV).$$
This sheaf is called the {\em group cohomology sheaf}.
We now consider the following composite of edge morphisms of the above spectral
sequences for $p \ge 1$
\begin{equation}\label{eqss}
H^p(\Gamma,H^0(U,\cV)) \to R^p(H^0(U,\cdot)^\Gamma)(\cV) \to H^0(U,\calH^p(\Gamma,\cV)).
\end{equation}
Also by \cite{Weibel}, Theorem~5.8.3, the edge morphisms are the natural maps.
If $W \subseteq U$ is an open set, then the restriction
$\res^U_W: \cV(U) \to \cV(W)$ induces natural maps on each of the three objects,
which are also denoted by $\res^U_W$.
It can be checked that the diagram
$$ \xymatrix@=1cm{
& H^p(\Gamma,H^0(U,\cV))        \ar@{->}[r] \ar@{->}[d]^{\res^U_W}  
& R^p(H^0(U,\cdot)^\Gamma)(\cV)  \ar@{->}[r] \ar@{->}[d]^{\res^U_W}
& H^0(U,\calH^p(\Gamma,\cV))                \ar@{->}[d]^{\res^U_W}\\
& H^p(\Gamma,H^0(W,\cV))        \ar@{->}[r] 
& R^p(H^0(W,\cdot)^\Gamma)(\cV)  \ar@{->}[r] 
& H^0(W,\calH^p(\Gamma,\cV))}$$
is commutative.
In other words, the edge morphisms in Eq.~(\ref{eqss}) give morphisms of presheaves
$$ \big(U \mapsto H^p(\Gamma,H^0(U,\cV))\big) \to  
   \big(U \mapsto R^p(H^0(U,\cdot)^\Gamma)(\cV)\big) \to \calH^p(\Gamma,\cV).$$
We call the first presheaf the {\em group cohomology presheaf}.
This terminology is justified, as the sheafification coincides
with the group cohomology sheaf in the cases of interest in the
present context.

\begin{proposition}\label{sheafi}
Assume that $\Gamma$ is of type $(\textrm{FP})_\infty$ over~$R$ 
(cf.~\cite{Bieri}, p.~6).
The morphisms of presheaves above become isomorphisms on the
sheafification.
In particular, for $x \in X$  one has
$\calH^i(\Gamma,\cV)_x \cong H^i(\Gamma,\cV_x)$ for all $i \in \NN$.
\end{proposition}

\pf
It follows directly from the definition that the second map is
sheafification. Indeed, let $\cV \hookrightarrow \cI^\bullet$
be an injective resolution.
Then 
$$\calH^p(\Gamma,\cV) = \Ker((\cI^p)^\Gamma \to (\cI^{p+1})^\Gamma)/
\Image((\cI^{p-1})^\Gamma \to (\cI^p)^\Gamma)$$
is by definition the sheafification of
$$U \mapsto R^p(H^0(U,\cdot)^\Gamma)(\cV) = 
\Ker(\cI^p(U)^\Gamma \to \cI^{p+1}(U)^\Gamma)/
\Image(\cI^{p-1}(U)^\Gamma \to \cI^p(U)^\Gamma).$$

As taking stalks is exact, we have for $x \in X$ that
$\cV_x \hookrightarrow \cI^\bullet_x$
is exact in the category of $R[\Gamma]$-modules.
We claim that this is a $\Gamma$-acyclic resolution of~$\cV_x$.
By Proposition~\ref{propss}~(\ref{pd}) we know that for all $U \subset X$
open $\cI^i(U)$ is an injective $R[\Gamma]$-module for all $i \ge 0$,
and hence that $H^q(\Gamma, \cI^i(U)) = 0$ for all $q \ge 1$.
Under the assumption by \cite{Bieri}, Proposition~2.4, we know that 
the functor $H^i(\Gamma,\cdot)$ commutes with direct limits,
whence, indeed,  
$$H^q(\Gamma,\cI^i_x) = \ilim_{U \ni x} H^q(\Gamma, \cI^i(U)) = 0$$
for all $q\ge 1$ and all $i\ge 0$, as claimed.
From this the particular statement follows directly, as the stalk at~$x$
in the center equals the cohomology of $\cV_x \hookrightarrow \cI^\bullet_x$,
which by the preceding computation coincides with 
$H^i(\Gamma, \cV_x) = \ilim_{U \ni x} H^i(\Gamma, H^0(U,\cV))$.
As isomorphism of sheaves can be tested on the stalks, the proposition follows.
\qed

\subsection*{A spectral sequence for Hecke triangle surfaces}

We again assume Notation~\ref{notation} and we let $\cC \in \{\HH, \Hbar\}$
and $X = \Gamma \backslash \cC \in \{Y_\Gamma, X_\Gamma\}$.
\smallskip

Let us recall some facts on~$\Hbar$ that we will use in the sequel.
The topology on $\Hbar$ extends the topology of~$\HH$ and
is obtained as follows. 
Using the action of~$G$ it suffices to give a system of open 
neighborhoods of the cusp~$\infty$, which is provided by the sets
$U_T = \{ a+ib \;| \; b > T \} \cup \{ \infty \}$ for every real $T > 0$.
Clearly, the intersection with any open set in $\HH$ 
is an open set in~$\HH$. 

The $G$-orbit of every $x \in \Hbar$ is a discrete set.
Around every $x \in \Hbar$ there even exists an open set~$U$
such that $gU \cap U \neq \emptyset$ implies that $gx = x$.
Indeed, this holds on $\HH$ (by the existence of a fundamental
domain) and we only need to check it on the cusps.
Let $\tau = x + i y$ with $y > 1$. By looking at the standard fundamental
domain we see that either $\textrm{Im}(g\tau) = \textrm{Im}(\tau)$,
implying $g\infty = \infty$, or $\textrm{Im}(g\tau) \le 1$.
Hence, with $T>1$ we have for any $g,h \in G$ that $h U_T \cap gU_T = \emptyset$,
unless $g\infty = h\infty$.

Next, we claim that $\Hbar$ is simply connected.
It suffices to show that any loop $L$ starting and ending
in $\infty$ is contractible. We may assume that it does not
pass through any other cusp (otherwise, we cut the loop into several
loops each one meeting only one cusp).
The following homotopy works
$$ g : [0,1] \times [0,1] \to \Hbar, (s,x) \mapsto 
\begin{cases}
 \frac{L(x)}{s} & \text{ if }s \neq 0 \text{ and } L(x) \neq \infty,\\
 \infty  & \text{ otherwise.}
\end{cases}$$
The map $g$ is continuous, as the preimage of $U_T$ is 
$[0,1]\times [0,1] \cap \{ (s,x) | \mathrm{Im}(L(x)) > sT \}$,
and, thus, is open in $[0,1]\times [0,1]$.
Clearly, $g(0,x)$ is the constant path from $\infty$ to $\infty$, whereas
$g(1,x) = L(x)$ is the loop we started with.
\smallskip

Denote by $\underline{V}$ the constant sheaf on~$\cC$ associated
with~$V$ together with its natural $\Gamma$-action, 
i.e.\ for an open set $U \subset \cC$ we let 
$\underline{V}(U) = \Hom_\cts(U,V)$ (equipping $V$ with the discrete
topology) together with isomorphisms $\phi_g: g_*\underline{V} \to \underline{V}$
for each $g \in \Gamma$ which on~$U$ are given by 
$$\Hom_\cts(g^{-1}U,V) \to \Hom_\cts(U,V), \;\; 
f \mapsto (u \mapsto g f( g^{-1}u) \; \forall u \in U).$$
We have that $\pi_*\underline{V}$ is in $\Shf_{X}(R[\Gamma]-\Mod)$.

\begin{lemma}\label{lempistar}
For any point $y \in \cC$ and any sheaf $\cF \in \Shf_\cC(R-\Mod)$
taking disjoint unions into products there is an isomorphism
$$ (\pi_*\cF)_{\pi(y)} \cong \prod_{\gamma \in \Gamma/\Gamma_y} 
\cF_{\gamma y}$$
of $R$-modules.
In particular, $\pi_*$ is an exact functor and for all $i \ge 0$
$$ H^i(\cC,\cF) \cong H^i(X, \pi_* \cF).$$
Moreover, there is an isomorphism
$$ (\pi_*\underline{V})_{\pi(y)}
\cong \Coind_{\Gamma_y}^\Gamma V$$
of $R[\Gamma]$-modules.
\end{lemma}

\pf
This follows from the fact that around any $x \in \cC$ there is an open
set~$U$ such that for any $\gamma \in \Gamma$ the intersection $\gamma U \cap U$
is empty, unless $\gamma x = x$. 
\qed

\begin{corollary}\label{rqg}
Let $\cF \in \Shf_X(R[\Gamma]-\Mod)$ and suppose $H^i(X,\cF) = 0$ for all $i \ge 1$.
Then 
$$H^q(\Gamma,H^0(X,\cF)) \cong R^q(H^0(X,\cdot)^\Gamma )(\cF)$$ 
for all $q \ge 0$ and, in particular, 
$$H^q(\Gamma,V) \cong R^q(H^0(X,\cdot)^\Gamma )(\pi_* \underline{V}).$$
\end{corollary}

\pf
The assumptions mean that the spectral sequence in Eq.~(\ref{spseins}) degenerates.
In the special case this is true since 
$H^i(X,\pi_*\underline{V})= H^i(\cC, \underline{V})= 0$ for all $i \ge 1$
by Lemma~\ref{lempistar} and the fact that $\cC$ is simply connected
and $\underline{V}$ is constant.
\qed

\begin{corollary}\label{sky}
The stalk in $x \in X$ of the group cohomology sheaf 
$\calH^q(\Gamma,\pi_*\underline{V})$ 
is $H^q(\Gamma_y,V)$ for any $y \in \cC$ with $\pi(y)=x$.
In particular, $\calH^q(\Gamma,\pi_*\underline{V})$
is a skyscraper sheaf on $X$ for all $q \ge 1$.
\end{corollary}

\pf
The group $\Gamma$ is of type $(\mathrm{FP})_\infty$. Indeed, by
\cite{Bieri}, Proposition~2.13 we know that free products of
groups of type $(\mathrm{FP})_\infty$ are of type $(\mathrm{FP})_\infty$.
As finite groups are clearly $(\mathrm{FP})_\infty$ (see \cite{Bieri},
Example~2.6), it follows that $G$ is. Finally, by \cite{Bieri}, Proposition~2.5,
subgroups of finite index in groups of type $(\mathrm{FP})_\infty$
are $(\mathrm{FP})_\infty$, whence $\Gamma$ is.

Thus, we may apply Proposition~\ref{sheafi}. The first statement
now follows from Lemma~\ref{lempistar} and Shapiro's lemma.
The special case is a consequence of the fact
that the non-trivially stabilized points of~$\cC$
for the action of $\Gamma$ are discrete.
\qed

Let us note as a consequence of the case $q=0$ that the sheaf 
$(\pi_*\underline{V})^\Gamma$ is locally constant on~$X$ 
if and only if $V^{\Gamma_y} = V$ for all $y \in \cC$. 

\begin{lemma}\label{edge}
The composition of the edge morphisms in Eq.~(\ref{eqss}) for $p \ge 1$
$$ H^p(\Gamma,V) \to R^p(H^0(X,\cdot)^\Gamma)(\pi_*\underline{V})
\to H^0(X, \calH^p(\Gamma,\pi_*\underline{V}))$$
is the restriction map from the theory of group cohomology, when we identify
$H^0(X, \calH^p(\Gamma,\pi_*\underline{V})$ with
$\prod_{x \in X} H^p(\Gamma_{y_x}, V)$
for a choice of $y_x \in \cC$ with $\pi(y_x)=x$.
\end{lemma}

\pf
A reformulation of Corollary~\ref{sky} is that
$\calH^p(\Gamma,\pi_*\underline{V})_x \cong H^p(\Gamma, \Coind_{\Gamma_{y_x}}^\Gamma V)$.
The composite edge morphism is given as follows. Choose injective resolutions
$(\pi_* \underline{V})(X) = V \hookrightarrow i^\bullet$ in $R[\Gamma]-\Mod$ and 
$\pi_* \underline{V} \hookrightarrow \cI^\bullet$ in $\Shf_X(R[\Gamma]-\Mod)$,
so that by the proof of Proposition~\ref{sheafi} and Lemma~\ref{lempistar}, 
$\Coind_{\Gamma_{y_x}}^\Gamma V = (\pi_* \underline{V})_x \hookrightarrow \cI^\bullet_{y_x}$
is a $\Gamma$-acyclic resolution.
Then the edge morphism followed by the identification in the statement
is induced from the restriction map on the sheaves
$(\pi_* \underline{V})(X) \to (\pi_* \underline{V})_x$,
which is the diagonal map $V \to \Coind_{\Gamma_{y_x}}^\Gamma V$.
On group cohomology the map in question is hence the natural homomorphism
$H^p(\Gamma,V) \to H^p(\Gamma,\Coind_{\Gamma_{y_x}}^\Gamma V)$.
Composing it with the isomorphism from Shapiro's lemma we obtain the
group theoretic restriction
$H^p(\Gamma,V) \to H^p(\Gamma_{y_x},V)$.
\qed

We have now established the following theorem.

\begin{theorem}\label{sps}
In Notation~\ref{notation} with $\cC \in \{\HH, \Hbar\}$
and $X = \Gamma \backslash \cC \in \{Y_\Gamma, X_\Gamma\}$, 
there is a spectral sequence
$$ H^p(X, \calH^q(\Gamma,\pi_*\underline{V})) \Rightarrow H^{p+q}(\Gamma,V),$$
in which the edge morphisms (for $p \ge 1$)
$$H^p(\Gamma,V) \to H^0(X, \calH^p(\Gamma, \pi_*\underline{V}))$$
become the group theoretic restriction map under the identification
$$H^0(X, \calH^q(\Gamma,\pi_*\underline{V}) \cong
\prod_{x \in X} H^q(\Gamma_{y_x}, V)$$
for a choice of $y_x \in \cC$ with $\pi(y_x)=x$.
\end{theorem}

\subsection*{Explicit description}

We keep assuming Notation~\ref{notation}.
Let $\cV$ be a sheaf of $R$-modules on~$Y_\Gamma$. 
From the Leray spectral sequence associated to~$j: Y_\Gamma \hookrightarrow X_\Gamma$ we
get the exact sequence
\begin{multline}\label{leray}
   0 \to H^1(X_\Gamma, j_* \cV) 
     \to H^1(Y_\Gamma,\cV)  
     \to H^0(X_\Gamma, R^1j_* \cV) \\
     \to H^2(X_\Gamma,j_* \cV) 
     \to H^2(Y_\Gamma, \cV).
\end{multline}
The {\em parabolic cohomology group (for $Y_\Gamma$ and $\cV$)} is
image of the map $\Hc^i (Y_\Gamma,\cV) \to H^i (Y_\Gamma,\cV)$. It is 
denoted by $\Hpar^i(Y_\Gamma,\cV)$.
Moreover, we call $H^0(X_\Gamma, R^1j_* \cV)$ the
{\em boundary cohomology group (for $Y_\Gamma$ and $\cV$)}. 

\begin{proposition}\label{parj}
There is a natural isomorphism of $R$-modules 
$\Hpar^1(Y_\Gamma,\cV) \cong H^1(X_\Gamma, j_* \cV).$
\end{proposition}

\pf
We consider the exact sequence of sheaves on~$X_\Gamma$
$$ 0 \to j_! \cV \to j_* \cV \to C \to 0,$$
in which $C$ is defined as the cokernel.
It is a skyscraper sheaf, as it is only supported
on the cusps. 
Hence, $H^1(X_\Gamma,C) = 0$ and the long exact sequence associated
to the short exact sequence of sheaves above yields that
the horizontal map is surjective in the commutative diagram
$$ \xymatrix@=.5cm{
  \Hc^1(Y_\Gamma,\cV) \ar@{->>}[r] \ar@{->}[rd] &
  H^1 (X_\Gamma,j_* \cV) \ar@{^{(}->}[d] \\
& H^1(Y_\Gamma,\cV),}$$
in which the vertical map comes from the Leray sequence Eq.~(\ref{leray}).
As it is injective, the proposition follows. 
\qed

We can now prove the principal result of this section.

\begin{theorem}\label{thmexpl}
Recall that we are assuming Notation~\ref{notation}.
Let $M$ denote the coinduced module $\Coind_\Gamma^{G} V$.
The following explicit descriptions hold:
$$ H^1(Y_\Gamma, (\pi_* \underline{V})^\Gamma) \cong 
M / \big( M^{\langle \sigma\rangle} + M^{\langle \tau \rangle} \big)$$
and
$$ \Hpar^1(Y_\Gamma, (\pi_*\underline{V})^\Gamma) \cong 
\ker\big(M / ( M^{\langle \sigma\rangle} + M^{\langle \tau \rangle}) 
\xrightarrow{1-\sigma} M/(1-T)M \big).$$
\end{theorem}

\pf
We first use that any non-trivially stabilized point $x$ of~$\Hbar$
is conjugate by some $g \in G$ to either $i$, $\zeta_n$ or~$\infty$.
The respective stabilizer groups are $G_i = \langle \sigma \rangle$,
$G_{\sigma \zeta_n} = \langle \tau \rangle$ and
$G_\infty = \langle T \rangle$.
Hence, from Mackey's formula (Corollary~\ref{corstabmackey})
we obtain
$$ \prod_{x \in Y_\Gamma} H^p(\Gamma_{y_x}, V) \cong
H^p(\langle \sigma \rangle, M) \oplus
H^p(\langle \tau \rangle, M)$$
and
$$ \prod_{x \in X_\Gamma} H^p(\Gamma_{y_x}, V) \cong
H^p(\langle \sigma \rangle, M) \oplus
H^p(\langle \tau \rangle, M) \oplus
H^p(\langle T \rangle, M),$$
where again $y_x \in \Hbar$ with $\pi(y_x) = x$.
From the spectral sequence from Theorem~\ref{sps} we now get
the exact sequences
$$ 0 \to H^1(Y_\Gamma, (\pi_* \underline{V})^\Gamma) \to H^1(\Gamma, V) 
  \to H^1(\langle \sigma \rangle, M) \oplus
H^1(\langle \tau \rangle, M)$$
and
$$ 0 \to H^1(X_\Gamma, (\pi_* \underline{V})^\Gamma) \to H^1(\Gamma, V) 
  \to H^1(\langle \sigma \rangle, M) \oplus
H^1(\langle \tau \rangle, M) \oplus
H^1(\langle T \rangle, M).$$
Using Proposition~\ref{parj} and comparing with the Mayer-Vietoris exact sequence 
(Proposition~\ref{mayervietoris}) yields the theorem.
\qed

\begin{remark}\label{RemMerel}
A study of the homology of modular curves (as Riemann surfaces) for a subgroup $\Gamma \le \PSL_2(\ZZ) =: G$
of finite index
was carried out by Merel in~\cite{MerelHecke}, also in order
to compute modular forms. 
His result \cite{MerelHecke}, Proposition~4,
$$H_1(X_\Gamma, \cusps, R) \cong M / \big( M^{\langle \sigma\rangle} + M^{\langle \tau \rangle} \big)$$
with $M = \Coind_\Gamma^{G} V$ and $\cusps = \Gamma \backslash G(\infty)$ 
is a special case of ours, as one can see as follows.

The general duality theorem (see \cite{Dold}, Proposition~VIII.7.2; note
that in this case \u{C}ech cohomology coincides
with singular cohomology (see e.g.\ \cite{Dold}, Proposition~VIII.6.12))
gives an isomorphism $H_1(X_\Gamma, \cusps, R) \cong H^1(Y_\Gamma, R)$,
so that we may apply Theorem~\ref{thmexpl}.
\end{remark}

\begin{remark}
The spectral sequence from Theorem~\ref{sps} has a geometric interpretation in
terms of analytic modular stacks (see \cite{Thesis}, II.3).
Due to the lack of suitable references for stacks we avoid their use here and
just state the result:

There is an analytic stack, denoted $[\Gamma \backslash \HH]$ together with
a projection in the category of analytic stacks $f: [\Gamma \backslash \HH] \to Y_\Gamma$.
The projection map $g: \HH \to [\Gamma \backslash \HH]$ allows one to define the
sheaf $(g_* \underline{V})^\Gamma$ on the analytic stacks, in analogy to the above
treatment. The derived functor cohomology for the functor taking global sections 
coincides with that of group cohomology, i.e.\
$$H^i([\Gamma\backslash \HH], (g_* \underline{V})^\Gamma) \cong H^i(\Gamma, V).$$
The spectral sequence from Theorem~\ref{sps} for $\cC = \HH$
can be identified with the Leray spectral
sequence associated to the projection map~$f$.
\end{remark}

\section{Comparison and computation of modular forms}\label{MFSec}

In this section we compare the various objects discussed so far and
sketch how they can be used for the computation of modular forms.

\subsection*{Comparison}

\begin{theorem}\label{compthm}
We assume Notation~\ref{notation}.
The following sequences are exact:
\begin{enumerate}
\item $0 \to H^1(Y_\Gamma, (\pi_* \underline{V})^\Gamma) \to H^1(\Gamma,V) \to
\prod_{x \in Y_\Gamma} H^1(\Gamma_{y_x},V) \to 0$,
\item $0 \to \Hpar^1(Y_\Gamma, (\pi_* \underline{V})^\Gamma) \to \Hpar^1(\Gamma,V) \to
\prod_{x \in Y_\Gamma} H^1(\Gamma_{y_x},V) \to 0$,
\item $M^G \to \prod_{x \in Y_\Gamma} \big(V^{\Gamma_{y_x}} / N_{\Gamma_{y_x}} V\big) \to \cM_R(\Gamma,V)
\to H^1(Y_\Gamma, (\pi_* \underline{V})^\Gamma) \to 0$,
\item $M^G \to \prod_{x \in Y_\Gamma}  \big(V^{\Gamma_{y_x}} / N_{\Gamma_{y_x}} V\big) \to \cCM_R(\Gamma,V)
\to \Hpar^1(Y_\Gamma, (\pi_* \underline{V})^\Gamma) \to 0$,
\end{enumerate}
where for all $x \in Y_\Gamma$ we have chosen $y_x \in \HH$ such that $\pi(y_x)=x$. 
As before, $M$ denotes $\Coind_\Gamma^G V$, which we identify with
$\Ind_\Gamma^G V$.
\end{theorem}

\pf
This follows from the explicit descriptions in Theorems~\ref{thmexpl}
and~\ref{ManinSymbols} and Propositions \ref{propcoh} and~\ref{propPara}.
\qed

\begin{corollary}\label{compcor}
We use the Notation~\ref{notation} with $R=\ZZ$.
The $\ZZ$-modules $H^1(\Gamma,V)$, $H^1(Y_\Gamma, (\pi_* \underline{V})^\Gamma)$
and $\cM_\ZZ(\Gamma,V)$ only differ by torsion.
The same statement holds for the $\ZZ$-modules
$\Hpar^1(\Gamma,V)$, $\Hpar^1(Y_\Gamma, (\pi_* \underline{V})^\Gamma)$
and $\cCM_\ZZ(\Gamma,V)$.
\end{corollary}

\begin{corollary}\label{compcoreins}
We assume Notation~\ref{notation} and
suppose that the order of $\Gamma_x$ is invertible in~$R$ for all $x \in \HH$.
Then there are isomorphisms
$$H^1(\Gamma,V) \cong H^1(Y_\Gamma, (\pi_* \underline{V})^\Gamma) \cong \cM_R(\Gamma,V)$$ 
and 
$$\Hpar^1(\Gamma,V) \cong 
\Hpar^1(Y_\Gamma, (\pi_* \underline{V})^\Gamma) \cong \cCM_R(\Gamma,V).$$ 
The statements hold, in particular, for the group $\Gamma_1(N)$ with 
$N \ge 4$.
\end{corollary}

\pf
This follows from Theorem~\ref{compthm}. For the last
part we use that under the condition $N \ge 4$ all $\Gamma_1(N)_x$ 
are trivial.
\qed

We should not fail to recall the following well-known facts.
The group $\Gamma_0(N)/\langle \pm 1 \rangle$ 
can only contain stabilizer groups of order $2$ or~$3$.
It contains no stabilizer of even order if and only if
$N$ is divisible by a prime~$q$ which is $3$ modulo~$4$ or by~$4$.
Furthermore, it does not contain any stabilizer of order~$3$
if and only if $N$ is divisible by a prime~$q$ 
which is $2$ modulo~$3$ or by~$9$ (for a proof see e.g.\ \cite{Thesis}).
Let us also mention that there are techniques for computing 
the torsion of the objects above
explicitly, see e.g.\ \cite{Thesis}, Proposition~2.4.8, and \cite{ArtZwei},
Proposition~2.6.

\subsection*{Computing modular forms}

We recall that a link between modular forms and the objects discussed in this
article is established by the Eichler-Shimura isomorphism
$$ \Hpar^1 (\Gamma,\CC[X,Y]_{k-2}) \cong S_k(\Gamma) \oplus 
\overline{S_k(\Gamma)},$$
where $k \ge 2$ is an integer,
$S_k(\Gamma)$ is the space of holomorphic weight~$k$ cusp forms for~$\Gamma$
and $\CC[X,Y]_{k-2}$ is the $\CC$-vector space of
homogeneous polynomials of degree~$k-2$.
This isomorphism exists for any Fuchsian group of the first kind
with parabolic elements (see \cite{Knopp}, Theorem~A, and the discussion there).

If the group in question is a congruence subgroup of $\SL_2(\ZZ)$,
there is a theory of Hecke operators on modular forms and on the objects
studied in this article. The Eichler-Shimura isomorphism
is compatible with the Hecke action (see \cite{DiamondIm}, Theorem~12.2.2).
As is well known, this can be used for computing Hecke algebras
and, hence, coefficients of modular forms by appealing to
the isomorphism
$\Hom_\ZZ (\TT, \CC) \xrightarrow{\phi \mapsto \sum_{n\ge 1} \phi(T_n) q^n} S_k(\Gamma)$
with the usual convention $q = e^{2 \pi i \tau}$
for $\tau \in \HH$.

For, the compatibility of the Hecke operators with
the Eichler-Shimura isomorphism implies that
the Hecke algebra of $S_k(\Gamma)$, i.e.\ the algebra generated by
the Hecke operators inside the endomorphism ring of $S_k(\Gamma)$,
is isomorphic to the Hecke algebra of 
$\Hpar^1(\Gamma,\CC[X,Y]_{k-2})$, and by Corollary~\ref{compcor},
also to the Hecke algebra of $\cCM_\QQ(\Gamma,\CC[X,Y]_{k-2})$ and
$\Hpar^1(Y_\Gamma,(\pi_* \underline{\CC[X,Y]_{k-2}})^\Gamma)$.
These Hecke algebras are finite dimensional, their dimensions
are known and the so-called Sturm bounds provide explicit bounds~$B$
such that the Hecke operators $T_1, T_2, \dots, T_B$ generate the
Hecke algebra (see \cite{SteinBook}).

It should be stressed that the torsion-free quotient of
$\Hpar^1(\Gamma,\ZZ[X,Y]_{k-2})$ is a natural $\ZZ$-structure in
$\Hpar^1(\Gamma,\CC[X,Y]_{k-2})$ (and similarly for the other objects
of this article). As the Hecke operators are already defined on this
$\ZZ$-structure, computations can be done in the integers. We also
obtain that the eigenvalues of the Hecke operators are algebraic
integers.

In the article \cite{ArtZwei} cases are described in which one may use
$\Hpar^1(\Gamma,\FF[X,Y]_{k-2})$ with a finite fields~$\FF$ for
computing Hecke algebras of Katz modular forms over~$\Fbar$.

\vspace*{.5cm}
\noindent Gabor Wiese\\
Universit\"at Duisburg-Essen\\
Institut f\"ur Experimentelle Mathematik\\
Ellernstr.~29\\
D-45326 Essen\\
Germany\\
E-mail: {\tt gabor@pratum.net}\\
Web page: {\tt http://maths.pratum.net}

\end{document}